\documentstyle[12pt]{article}
\begin{document}
\title{{On extensions of Myers' theorem}}
\author{Xue-Mei Li\thanks{Research supported  in part by NATO Collaborative
Research Grants Programme 0232/87 and by SERC grant GR/H67263.
1991 Mathematical subject classification 60H30,53C21}}
\date{}
\maketitle

\newcommand{\A}{{\bf \cal A}}
\newcommand{\B}{{ \bf \cal B }}
\newcommand{\C}{{\cal C}}
\newcommand{\F}{{\cal F}}
\newcommand{\G}{{\cal G}}
\newcommand{\h}{{\cal H}}
\newcommand{\K}{{\cal K}}
\newcommand{\half}{{  {1\over 2}  }}
\newcommand{\heatsemif}{{ {\rm e}^{ \half t\Delta ^{h,1}}   }}
\newcommand{\heatsemi}{{ {\rm e}^{\half t \Delta ^{h}}   }}

\newtheorem{theorem}{Theorem}
\newtheorem{proposition}[theorem]{Proposition}
\newtheorem{lemma}[theorem]{Lemma}
\newtheorem{corollary}[theorem]{Corollary}
\newtheorem{definition}{Definition}[section]

\def\limsup{\mathop{\overline{\rm lim}}}
\def\liminf{\mathop{\underline{\rm lim}}}
\def\exp{{\rm e}}

\begin{abstract}
Let $M$ be a  compact Riemannian manifold  and $h$ a smooth function on $M$.
Let $\rho^h(x)=\inf_{|v|=1}\left(Ric_x(v,v)-2Hess(h)_x(v,v) \right)$.  Here
 $Ric_x$ denotes the Ricci curvature at $x$ and $Hess(h)$ is the Hessian of $h$. 
 Then $M$ has finite fundamental group if $\Delta ^h-\rho^h<0$.
Here $\Delta ^h=: \Delta  +2L_{\nabla h}$ is the Bismut-Witten Laplacian.
This leads to a quick proof of recent results on extension of Myers'
theorem to manifolds with mostly positive curvature. There is also a similar
result for noncompact manifolds.

\end{abstract}

An early result of Myers says a complete Riemannian manifold with Ricci
 curvature bounded below by a positive number is compact and has finite
fundamental group. See e.g. \cite{Gallot}.
 Since then efforts have been made to get the same
type of result but to allow a little bit of negativity of the curvature
(see  B\'erard and Besson\cite{be-be}).
Wu \cite{WUJ}  showed that Myers' theorem
holds if the manifold is allowed to have negative curvature on a set of
small diameter, while Elworthy and Rosenberg \cite{EL-RO91}  considered
 manifolds with some negative   curvature on a set of small volume,
followed by  recent work of Rosenberg and Yang \cite{RO-YA}.  We use
a method of Bakry \cite{BA86} to obtain  a 
result given in terms of the potential kernel related to 
$\rho(x)=\inf_{|v|=1} Ric_x(v,v)$, which gives a quick probabilistic proof
 of recent results on extensions of Myers' theorem. Here $Ric_x$ denotes  the
Ricci curvature at $x$.

\vskip 20pt

Let $M$ be a complete Riemannian manifold, and $h$ a smooth real-valued
 function on it. Assume Ric$-2$Hess(h) is bounded from below, where Hess(h)
 is the hessian of $h$. Denote by $\Delta^h$ the Bismut-Witten Laplacian
 (with probabilistic sign convention) defined by:
$\Delta  h=\Delta+2 L_{\nabla h}$ on $C_K^\infty$ the space of smooth 
differential forms with compact support.  Here $L_{\nabla h}$ is the Lie 
derivative in direction of $\nabla h$. Then  the closure of $\Delta^h$ is a
 negative-definite self-adjoint differential operator on $L^2$ functions
 (or $L^2$ differential forms) with respect to $e^{2h}dx$ for $dx$ the 
standard Lebesgue measure on $M$. We shall use the same notation for
 $\Delta^h$ and its closure.   By the  spectral theorem there is a heat
 semigroup $P_t^h$ satisfying the  following heat equation:

$${\partial u(x,t)\over \partial t}=\half \Delta^hu(x,t).$$
We shall denote by $P_t^h\phi$ the solution with initial value $\phi$.
For clarity, we also use $P_t^{h,1}$ for the corresponding heat semigroup
for one forms. Then for  a function $f$ in $C_K^\infty$,

\begin{equation}
dP_t^hf=P_t^{h,1}(df).
\label{eq: 1}
\end{equation}

On $M$ there is a h-Brownian motion $\{F_t(x):t\ge 0\}$, i.e. a path continuous
strong Markov process with generator $\half \Delta ^h$ for each starting point
 $x$. For a fixed point $x_0\in M$, we shall write $x_t=F_t(x_0)$. 
 Then $P_t^hf(x)=Ef(F_t(x))$ for all bounded $L^2$ functions.

 Let $\{W_t^h(-), t\ge 0\}$ be the solution flow to the following covariant 
equation along  h-Brownian paths $\{x_t\}$:

\begin{equation}
\left\{ \begin{array}{ll}

{D\over \partial t}W_t^h(v_0) &=-\half Ric_{x_t}\left(W_t^h(v_0),-\right)^\#
         		+Hess(h)_{x_t}\left (W_t^h(v_0),-\right)^\#,\\
W_0^h(v_0)&=v_0, \hskip 24pt v_0\in T_{x_0}M.
\end{array}\right.
\label{eq: 2}
\end{equation}
Here $\#$ stands for the adjoint. The solution flow $W_t^h$ is called the
Hessian flow.

Let $\phi$ be a bounded 1-form, then for $x_0\in M$, and $v_0\in T_{x_0}M$
\begin{equation}
E\phi(W_t^h(v_0))=P_t^{h,1}\phi(v_0).
\label{eq: 3}
\end{equation}
if Ricci$-2$ Hess(h) is bounded from below. See  e.g. \cite{ELflour} and
 \cite{ELflow}.

Formula $(\ref{eq: 3})$ gives the following estimates on the heat semigroup:
\begin{equation}
|P_t^{h,1}\phi|\le |\phi|_{\infty} E|W_t^h|.
\label{eq: 4}
\end{equation}

Let $\rho^h(x)=\inf_{|v|= 1}\left(Ric_x(v,v)-2Hess(h)(x)(v,v)\right)$ and 
write $\rho$ for $\rho^h$ if $h=0$.
Then covariant equation (\ref{eq: 2}) gives:
\begin{equation}
E|W_t^h|\le E\exp^{-\half\int_0^t\rho^h(x_s)ds}
\label{eq: 5}
\end{equation}
as in \cite{ELflour}.

 Let $P_t^{\rho^h}$  be the $L^2$ semigroup generated by the  Schr\"odinger
operator  $\half \Delta ^h -\half \rho^h$.
 Then 
$$P_t^{\rho^h} f(x_0)=E\left[f(x_t)\exp^{-\half\int_0^t\rho^h(x_s)ds}\right]$$
by  the Feyman-Kac formula. So equation  (\ref{eq: 5}) is equivalent to
$$ E|W_t^h|\le P_t^{\rho^h} 1.$$

Let $Uf$ be the corresponding potential kernel defined by:
$$Uf(x_0)=\int_0^\infty E\left[f(x_t)\exp^{-\half\int_0^t\rho^h(x_s)ds}\right] dt.$$

Following Bakry's paper \cite{BA86}, we have the following theorem:

\begin{theorem}\label{th:1}
Let $M$ be a complete Riemannian manifold with Ric$-2$Hess(h) bounded from
 below. Suppose 
\begin{equation}
\sup_{x\in K} U1(x)<\infty
\label{eq: 6}
\end{equation}
for each compact set $K$. Then $M$ has finite h-volume(i.e.
$\int_M \exp^{2h(x)} dx<\infty$),  and finite fundamental group.
\end{theorem}

\noindent
{\bf Proof:} We follow \cite{BA86}.
Let $f\in C_K^\infty$, then $Hf=\lim_{t\to \infty}P_tf$ is an $L^2$ harmonic
function. Assume h-vol$(M)=\infty$, then $Hf=0$. We shall prove this is 
impossible.
Let $f$, $g$ $\in C_K^\infty$, then:

 \begin{eqnarray*}
&&\int_M(P_t^hf-f)g\exp^{2h} dx\\
&&=\int_M \int_0^t\left( {\partial \over \partial s} P_s^hf\right) g
\exp^{2h} dx\\
&&=\int_0^t \int_M <\nabla P_s^hf, \nabla g>\exp^{2h} dx ds\\
&&\le |\nabla f|_{\infty} \int_M |\nabla g| \left(\int_0^t E|W_s^h|ds\right)
\exp^{2h} dx,   \\
&&\le |\nabla f|_{\infty} \left(\sup_{x\in sup(g)} \int_0^\infty E|W_s^h|ds\right)
 |\nabla g|_{L^1} \\
&&\le c|\nabla f|_\infty |\nabla g|_{L^1}.
\end{eqnarray*}
Here $c=\sup_{x\in sup(g)}\left[U1(x)\right]
=\sup_{x\in sup(g)} 
\left(\int_0^\infty E\left(\exp^{-\int_0^t \rho^h(F_s(x))ds}\right)dt\right)$,
and $sup(g)$ denotes the support of $g$.

Next take $f=h_n$, for $h_n$ an increasing sequence of smooth functions 
approximating $1$ with $0\le h_n\le 1$ and $|\nabla h_n|\le {1\over n}$, see
e.g. \cite{BA86}.

Then 
$$\int_M(P_t^h h_n-h_n)g\exp^{2h} dx\le c{1\over n}|\nabla g|_{L^1}.$$
First let $t$ go to infinity, then let $n\to \infty$ to obtain:

$$-\int_M g\exp^{2h} dx\le 0.$$
This gives a contradiction with a suitable choice of $g$. So we conclude
 h-vol$(M)<\infty$.

Let  $p\colon \tilde M \to M$ be the universal covering space for $M$ with 
induced Riemannian metric on $\tilde M$.  For $p(\tilde x)=x$, let
 $\{\tilde F_t(\tilde x), t\ge 0\}$ be the horizontal lift of $\{F_t(x)\}$
to $\tilde M$.  Denote by $\tilde Ric$  the Ricci curvature on $\tilde M$,
 $~\tilde h$ the lift of $h$ to $\tilde M$, and $\tilde {\rho^h}$ the
 corresponding lower bound for $\tilde Ric-2Hess(\tilde h)$. Then the induced
  $\{\tilde F_t(\tilde x), t\ge 0\}$ is a h-Brownian motion on $~\tilde M$.
See e.g. \cite{ELflour}. Note also  $\tilde {\rho^h}$ satisfies
$$\sup_{\tilde x\in\tilde  K} \int_0^\infty 
E\left(\exp^{-\half\int_0^t\tilde {\rho^h}(\tilde F_s(\tilde x))ds}\right)dt
=\sup_{x\in p( \tilde K)} \int_0^\infty 
E\left(\exp^{-\half \int_0^t \rho^h( F_s(x))ds}\right)dt$$
  for any $\tilde K\subset \tilde M$ compact. The same calculation as above 
will show that $\tilde M$ has finite h-volume, therefore
 $p$ is a finite covering and so $M$ has  finite fundamental group.
\hfill \rule{3mm}{3mm}

\bigskip
Let $r(x)$ be the Riemannian distance between $x$ and a fixed point of $M$
 and take $h$ to be identically zero:
\begin{corollary}
Let $M$ be a complete Riemannian manifold with
\begin{equation}\label{eq: volume}
Ric_x>-{n\over n-1}{1\over r^2(x)}, \hskip 6pt \hbox{when} \hskip 4pt r>r_0
\end{equation}
for some $r_0>0$.
 Then the manifold is compact if
$\sup_{x\in K}U1(x)<\infty$ for each compact set $K$.
\end{corollary}

\noindent {\bf Proof:}
This is a consequence of the  result of \cite{ch-gr-ta}:
 A complete Riemannian  manifold with (\ref{eq: volume}) has infinite
 volume.  See also  \cite{Wu91}. \hfill\rule{3mm}{3mm}

For another extension of Myers' compactness theorem, see \cite{ch-gr-ta}
where a diameter estimate is also obtained.

\bigskip

In the following we shall assume $M$ is compact  and  get the following corollary:

\begin{corollary}\label{corollary}
Let $M$ be a compact Riemannian manifold and $h$ a smooth function on it.
 Then $M$  has  finite fundamental group if $\Delta ^h-\rho^h<0$.
\end{corollary}

\noindent
{\bf Proof:}
 Let $\lambda_0$ be the  minimal eigenvalue  of $\Delta ^h -\rho^h$.
Then 
$$\limsup_{t\to \infty} {1\over t} \sup_M \log 
E\left(\exp^{-\half \int_0^t\rho^h(x_s)ds}\right)\le \lambda_0<0.$$
See e.g. \cite{EL-RO88}.
Thus there is a number $T_0$ such that if $t\ge T_0$, 
$$\sup_M E\left(\exp^{-\half \int_0^t \rho^h(x_s)ds}\right)
 \le \exp^{\lambda_0 t}.$$
Therefore
\begin{eqnarray*}
&&\sup_M\int_0^\infty E\left(\exp^{-\half \int_0^t \rho^h(x_s)ds}\right)dt\\
&&\le \int_0^{T_0}
 E\left(\exp^{-\half \inf_{x\in M} \left[\rho^h(x)\right]t}\right)dt
+\int_{T_0}^\infty E\left(\exp^{-\half \int_0^t \rho^h(x_s)ds}\right)dt\\
&&<\infty.
\end{eqnarray*}
The result follows from theorem \ref{th:1}.
\hfill\rule{3mm}{3mm}

Let $d_h=\exp^{h}d\exp^{-h}$. It has adjoint 
$\delta_h=\exp^{-h}\delta\exp^{h}$ on
$L^2(M,dx)$. Let $\Box_h$ be the Witten Laplacian defined by:
$$\Box_h=-(d_h+\delta_h)^2.$$
By conjugacy of $\Box_h$ on $L^2(M,dx)$ with $\Delta ^h$ on
 $L^2(M,e^{2h}dx)$,
the condition "$\Delta ^h-\rho^h<0$" becomes:
$$\Box_h-\rho^h<0$$
on $L^2(M,dx)$. On the other hand,
$$\Box_h=\Delta -||dh||^2-\Delta  h.$$
See e.g.  \cite{EL-RO92}. This gives: a compact manifold 
has finite fundamental  group if 
$$\Delta -||dh||^2-\Delta  h-\rho^h<0$$
on $L^2(M,dx)$.

\bigskip

Corollary \ref{corollary} leads to the following theorem from  
\cite{RO-YA}: 
Let ${\cal N}={\cal N}(K,D,V,n)$ be the collection of n-dimensional
Riemannian manifolds with
Ricci curvature bounded below by $K$, diameter bounded above by $D$, and volume
bounded below by $V$.

\begin{corollary}[Rosenberg\& Yang]
Choose $R_0>0$. There exists $a=a({\cal N}, R_0)$ such that a manifold
 $M\in {\cal N}$ with ${\rm vol}\{x: \rho(x)<R_0\} <a $
 has finite fundamental group. Here "vol" denotes the volume of the relevant 
set.
\end{corollary}

\noindent
{\bf Proof:} Let $h=0$ in corollary \ref{corollary}. Then under the
 assumptions in the corollary, $\Delta -\rho<0$ according to
 \cite{EL-RO91}.  The conclusion follows from
corollary \ref{corollary}.
\hfill\rule{3mm}{3mm}

\vskip 24pt

\noindent
{\bf Acknowledgement:} The author is grateful to Professor D. Elworthy and 
Professor S. Rosenberg
for helpful comments and encouragement.

\noindent
Address:	

\noindent
Mathematics Institute,  University of Warwick,

\noindent
Coventry CV4,7AL

\noindent U.K.

\noindent
e-mail: xl@maths.warwick.ac.uk

\end{document}